\documentclass[11pt]{amsart}
\newtheorem{theorem}{Theorem}[section]
\newtheorem{lemma}[theorem]{Lemma}
\newtheorem{proposition}[theorem]{Proposition}
\newtheorem{corollary}[theorem]{Corollary}

\theoremstyle{definition}
\newtheorem{definition}[theorem]{Definition}
\newtheorem{example}[theorem]{Example}

\theoremstyle{remark}
\newtheorem{remark}[theorem]{Remark}

\numberwithin{equation}{section}

\begin{document}

\title
{Quaternionic bundles and betti numbers of Symplectic 4-manifolds
with Kodaira dimension zero}

\author{Tian-Jun Li}
\address{School  of Mathematics\\  University of Minnesota
\\ Minneapolis, MN 55455\\}
\email{tjli@math.umn.edu}
\thanks{The author is supported in part by NSF grant 0435099 and the McKnight fellowship.}

\maketitle

\section{Introduction}

For a minimal symplectic $4-$manifold $M$ with symplectic form $\omega$
and
symplectic canonical class $K_{\omega}$, the Kodaira dimension of $(M,\omega)$
 is defined in the following way (\cite{L}, \cite{McS}):

\[ \kappa(M,\omega)=\left\{ \begin{array}{ll}
-\infty &\hbox{if $K_{\omega}\cdot [\omega]<0$ or $K_{\omega}\cdot K_{\omega}<0$},\\
0&\hbox{if $K_{\omega}\cdot [\omega]=0$ and $K_{\omega}\cdot K_{\omega}=0$},\\
1&\hbox{if $K_{\omega}\cdot [\omega]> 0$ and $K_{\omega}\cdot K_{\omega}=0$},\\
2&\hbox{if $K_{\omega}\cdot [\omega]>0$ and $K_{\omega}\cdot K_{\omega}>0$}.\\
\end{array}
\right. \]
 The Kodaira dimension of a non-minimal manifold is
defined to be that of any of its minimal models.

It is shown in \cite{L} that, if $\omega$ is a K\"ahler form on a
complex surface $(M, J)$, then $\kappa(M, \omega)$ agrees with the
usual holomorphic Kodaira dimension of $(M, J)$.

It is also shown in \cite{L}  that minimal symplectic $4-$manifolds
with $\kappa=0$ are exactly those with torsion canonical class, thus
can be viewed as {\it symplectic Calabi-Yau surfaces}. Known
examples of  symplectic $4-$manifolds with torsion canonical class
are either K\"ahler surfaces with (holomorphic) Kodaira dimension
zero or $T^2-$bundles over $T^2$ (\cite{Ge}, \cite{L}). They all
have small Betti numbers and Euler numbers:
 $b^+\leq 3$, $b^-\leq 19$ and $b_1\leq 4$; and the Euler number is between
$0$ and $24$. It is speculated in \cite{L} that these are the only ones.
In this paper we prove that it is true up to rational homology.

\begin{theorem}  \label{main} Suppose $M$ is a {\it minimal} symplectic $4-$manifold with $\kappa=0$. Then the rational homology of $M$ is the
same as that of K3 surface, Enriques surface or a $T^2-$bundle over
$T^2$. In particular, we have

\begin{enumerate}

\item the Euler number of $M$ is $0, 12$ or $24$,

\item the signature of
$M$ is $-16, -8$ or $0$,  and

\item  the Betti numbers of $M$ satisfy the following bounds:
$$b^+\leq 3, \quad b^-\leq 19, \quad b_1\leq 4.$$

\end{enumerate}

\end{theorem}

The case $b_1=0$ follows from \cite{MS}.
Under the additional assumption that $b_1\leq 4$, this theorem is proved in \cite{L}.
The key is really to bound $b^+$.
Our approach here is similar to that  in \cite{L}, which  is to
show that, on a closed smooth oriented $4-$manifold with $2\chi+3\sigma=0$ and $b^+>3$,
 the mod 2 Seiberg-Witten
invariant of any reducible Spin$^c$ structure vanishes. In this
paper we will call a Spin$^c$ structure reducible if it admits a
reduction to a spin structure. We have mentioned that {\it minimal}
symplectic $4-$manifolds with Kodaira dimension zero are exactly
those with torsion symplectic canonical class. In addition, a closed
symplectic $4-$manifold with $b^+>1$ and torsion canonical class
actually has trivial canonical class, and hence is a spin manifold.
For spin manifolds there are  stable cohomotopy and stable
homotopy/framed bordism refinements of the Seiberg-Witten invariants
of spin manifolds in \cite{BF}, \cite{F} and \cite{FL}, which take into account of
the $Pin(2)$ symmetry of the Seiberg-Witten equations.
Such refinements are used in section \ref {homotopy and bordism} to
construct an unoriented bordism SW invariant when $b^+\geq 2$
following \cite{FKMM}. The main theorem follows from a rather
general vanishing result of the unoriented bordism SW invariant. The
proof of the vanishing result relies on a few properties of
quaternionic bundles proved in section \ref {H}, which certainly are
of independent interest.

A basic conjecture of Gompf in \cite{Go} is that a symplectic
$4-$manifold with $\kappa$ at least zero has non-negative Euler
number. Theorem \ref{main} confirms it when  $\kappa=0$.

\begin{corollary}\label{nonminimal}
Any symplectic $4-$manifold with $\kappa=0$ has $b^+\leq 3, b_1\leq
4$ and non-negative Euler number.
\end{corollary}

Notice that the bound for $b_1$ is the same as the dimension. One
could speculate whether such a bound continues to hold in higher
dimensions.


Part of this work was completed during the visits at IPAM in March
2003, University of Tokyo in December 2004 and June 2005. The author
would like to thank these institutions for hospitality. The author
is particular grateful to  M. Furuta for many helpful discussions.
We also appreciate B. Gompf, B. H. Li, Y. Ruan, P. Seidel, B.
Siebert, A. Stipsicz and S. T. Yau's interest in this work. Finally
we would like to thank a referee for the many useful suggestions.
After we completed this manuscript we learned the preprint \cite{B2}
where the Betti number bounds are also proved. This research is
partially supported by NSF and the McKnight fellowship.

\section{Quaternionic bundles} \label{H}

Let $J$ be a smooth manifold with an involution $\iota_J$ and with
nonempty and isolated fixed point set.

\begin{example}We are interested in the case that $J$ is the torus
$T^m={\mathbb R}^m/{\mathbb Z}^m$ with  $\iota_J$  given by $x\to
-x$ using the coordinates of ${\mathbb R}^m$. In this case we use
$O_J$ to denote the image of the origin in ${\mathbb R}^m$. Notice
that there are $2^m$ fixed points including $O_J$.
\end{example}

Recall that a bundle map between complex bundles is called
anti-complex if it anti-commutes with the multiplication by ${\bf
i}=\sqrt {-1}$.

\begin{definition} A complex bundle $Q$ on $J$ with an anti-complex
lift $\iota_Q$ of $\iota_J$ is called a {\it quaternionic} vector
bundle if $\iota_Q\circ \iota_Q=-1$.
\end{definition}

 Since the fixed point set
is nonempty and the fiber over any fixed point is a  space over the
quaternions ${\mathbb H}={\mathbb C}\oplus{\mathbb C}{\bf j}$, the
rank is necessarily even. However, we should warn the readers that a
{\it quaternionic} bundle here is not a bundle over ${\mathbb H}$.
In particular, the rank of a {\it quaternionic} bundle is its rank
as a complex bundle. The Grothendick group of the {\it quaternionic}
vector bundles is denoted by $KQ(J)$ (first appeared in \cite{D}).

Let $\underline {\mathbb H}^l$ be the rank $2l$ quaternionic vector
bundle  $J\times {\mathbb H}^l$ with the anti-complex map
$\iota_{\underline {\mathbb H}^l}:(x, q)\longrightarrow (\iota_J x,
q{\bf j})$, where $q{\bf j}$ is the right multiplication of $q$ by
the quaternion number ${\bf j}$. A rank $2l$ quaternionic vector
bundle $E$ is called trivial if there is a complex isomorphism
$\Phi:E\to \underline{\mathbb H}^l$ with $\Phi\circ
\iota_E=\iota_{\underline {\mathbb H}^l}\circ \Phi$. Quaternionic
vector bundles over low dimensional tori are classified in
\cite{FKMM}, and $KQ(T^m)$ is calculated in \cite{FK}.

Just as complex vector bundles are
acted upon by $U(1)$ via the complex multiplication,
quaternionic vector bundles are  naturally acted upon
by the group $Pin(2)$,  which is
 generated by $U(1)$ and the symbol $\iota$ with
the relations
$$\iota^2=-1,\qquad \iota z\iota^{-1}=z^{-1} \quad\hbox{for} \quad z\in U(1). $$
Clearly $Pin(2)$ fits into  the short exact sequence
$$1\longrightarrow U(1)\longrightarrow Pin(2)\longrightarrow \{\pm
1\}\longrightarrow 1.$$ Notice that  $Pin(2)$ is isomorphic to the
subgroup of ${\mathbb H}$ generated by $U(1)=\{\cos\theta+ {\bf
i}\sin\theta\}$ and ${\bf j}$.

We first specify the $Pin(2)$ action  on $J$: it is simply defined
 via the surjection of $Pin(2)$  onto the order $2$ group $\{id, \iota_J\}$.
For a quaternionic vector bundle $E$ over $J$, since $\iota_E$ is anti-complex,
 $Pin(2)$ acts on $E$ via the complex multiplication and $\iota_E$.

\begin{remark} \label{real} We will also need the simple fact that, for a real vector space $W$, the  trivial real vector bundle
$\underline W=J\times W$ is also $Pin(2)-$equivariant via the
involution $\iota_W:(x, a)\longrightarrow (\iota_J x, -a)$ and the
surjection $Pin(2)\longrightarrow \{id, \iota_W\}$.
\end{remark}

Notice that, since $Pin(2)$ is compact, there exists a
$Pin(2)-$invariant Hermitian inner product on any
$Pin(2)-$equivariant bundle.

The two main results about quaternionic bundles are Theorems
\ref{split} and \ref{monomorphism}. The first one is about splitting
off a trivial summand. We start with the following characterization.

\begin{lemma} \label{section} Let $E$ be a quaternionic bundle. Then
$E$ splits into $\underline {\mathbb H} \oplus E'$ for some
quaternionic bundle $E'$ if and only if there is a nowhere vanishing
section $s$ such that $s$ and $\iota_E s$ are complex linearly
independent everywhere.

\end{lemma}

\begin{proof} Suppose $E$ splits into $\underline {\mathbb H} \oplus E'$
for some  quaternionic bundle $E'$.  The constant section
${\underline 1}=(x,1)$ of $\underline {\mathbb H}$ is a nowhere
vanishing section of $E$, which we call $s$. Notice that the
constant section ${\underline {\bf j}}=(x,{\bf j})$ can be also
written as $\iota_{\underline {\mathbb H}} { \underline 1}$ .
Therefore, due to the $Pin(2)-$equivariance, the section $\iota_E s$
corresponds to ${\underline {\bf j}}$. Since $1$ and ${\bf j}$ form
a complex basis of $\mathbb H$, $s$ and $\iota_E s$ are complex
linearly independent at every point.

Conversely, we obtain a map from $s$ a quaternionic map
$$\phi:\underline
{\mathbb H}\to E, \quad  (x, a+b{\bf j})\to  as_x+b(\iota_Es)_x,$$
where $a, b\in {\mathbb C}$.  $\phi$ is an embedding because  $s_x$
and $(\iota_E s)_x$ are complex linearly independent for any $x$.
The required splitting is then obtained by choosing a
$Pin(2)-$invariant Hermitian metric and letting $E'$ be the
orthogonal complement of $\phi({\mathbb H})$.
\end{proof}

It is certainly not true that if $E$ has a nowhere vanishing
section,
 then it has one such section $s$ such that
 $s$ and $\iota_E s$ are
complex linearly independent everywhere. Otherwise, since  every
rank 2 quaternionic bundle over $T^2$ has a nowhere vanishing
section by dimension reason, we would draw the conclusion that every
such bundle is isomorphic to $\underline {\mathbb H}$. But by the
classification of quaternionic bundles over low dimensional tori
\footnote{up to dimension $4$} in \cite{FKMM}, there is a (unique)
non-trivial rank 2 quaternionic bundle over $T^2$.

For a nowhere vanishing section $s$, clearly $s$ and $\iota_Es$ are
complex linearly independent over any fixed point of $J$. On the
other hand, if $x$ is not a fixed point of $J$, then $s$ and
$\iota_Es$ are complex linearly independent over $x$ if and only if
$s_x$ is not mapped by $\iota_E$ to a point in the complex line
generated by $s_{\iota_Jx}$.
 To further investigate this problem for a quaternionic bundle of
rank $2l$ we introduce some auxiliary bundles.

{\bf The (complex) projective space bundle $P(E)$}. Let $P(E)$
denote the (complex) projective space bundle of $E$, which is a
${\mathbb C}P^{2l-1}-$bundle over $J$.  For any nonzero $u\in E_x$,
we use $[u]\in P(E)|_x$ to denote the complex line generated by $u$.
Suppose $s$ is a nowhere vanishing section of $E$, then it defines a
section $[s]$ of $P(E)$. Notice that $\iota_E$ sends a complex line
in $E$ to a complex line and therefore induces an action on $P(E)$.
This is simply because, for any nonzero $u\in E_x$, we have
$$\iota_E ((a+bi)u)=(a-bi)\iota_E(u).$$
In other words,  $\iota_E[u]=[\iota_Eu]$, where we continue to use $\iota_E$ to
denote the induced action on $P(E)$. Clearly $\iota_E$ is an
involution on $P(E)$.

{\bf The quaternionic bundle $\hat E$}. Let $\hat E$ be the pull
back bundle of $E$ under $\iota_J$. The fiber of $\hat E$ over $x$
is the fiber of $E$ at $\iota_Jx$, and vice versa. we can define the
quaternionic structure on $\hat E$ by
 $$\iota_{\hat E} (x,
v)=(\iota_Jx, \iota_E|_{\iota_Jx}(v)),$$ although we will not use
this structure on $\hat E$.  For a section $s$ of $E$, we have the
associated section $\hat s$ of $\hat E$, which is defined to be
$$\hat s_x=s_{\iota_Jx}$$ for any $x$.

{\bf The bundle  $E\oplus {\hat E}$ and the involution $\tau$}.
Consider the direct sum quaternionic bundle $E\oplus {\hat E}$.
 Since
$$E|_x={\hat E}|_{\iota_Jx}, \quad {\hat E}|_x=E|_{\iota_Jx},$$
$E\oplus {\hat E}$ has  an involution
$$\tau:(x, u, v)\longrightarrow  (\iota_Jx,  v,  u), \quad u\in E_x, v\in {\hat E}_x,$$
which covers the involution $\iota_J$ of $J$.
 Observe that  a
$\tau-$invariant section of $E\oplus {\hat E}$ is nothing but a
${\mathbb Z}_2-$equivariant map from $J$ to $E\oplus {\hat E}$.
Notice that  sections of $E$ correspond exactly to $\tau-$invariant
sections of $E\oplus \hat E$. On the one hand, any section $s$ of
$E$  gives rise to a $\tau-$invariant section $(s, \hat s)$ of
$E\oplus {\hat E}$. On the other hand, if a section  $(f, g)$ of
$E\oplus {\hat E}$ is $\tau-$invariant, then, for any $x\in J$, we
have
$$(\iota_Jx, f_{\iota_Jx}, g_{\iota_Jx})=\tau(x, f_{x},
g_{x})=(\iota_Jx, g_x, f_x).$$ Thus $g_x=f_{\iota_Jx}$, and in
particular, $g$ is completely determined by $f$.

{\bf The fiber product $P(E)\times_J P({\hat E})$}.  Consider the
fiber product $P(E)\times_J P({\hat E})$, which is a bundle with
fiber ${\mathbb C}P^{2l-1}\times {\mathbb C}P^{2l-1}$. The
involution $\tau$ on $E\oplus \hat E$ also induces an involution on
$P(E)\times_J P({\hat E})$, still denoted by $\tau$. A fixed point
of $\tau$ is of the form $(x, [u], [u])$, where $x$ is a fixed point
of $\iota_J$. A nowhere vanishing section $s$ of $E$ gives rise to a
$\tau-$invariant section $([s], [\hat s])$ of $P(E)\times_J P({\hat
E})$,  which can be viewed as  an ${\mathbb Z}_2-$equivariant map
from $J$ to  $P(E)\times_J P({\hat E})$.


{\bf The submanifold $D\subset P(E)\times_J P({\hat E})$}. Consider
the subset of $P(E)\times_J P({\hat E})$,
$$D=\{(x, [u], [v])\in P(E)\oplus P({\hat E})|[v]=\iota_E [u]\}.$$
Notice that if $(x, [u], [v])\in D$, then $[u]=\iota_E [v]$ as well.
This implies that $D$ is $\tau-$invariant. However $D$ does not
contain any fixed points of $\tau$.
Observe also that $D$
is diffeomorphic to $P(E)$ via the map $$(x, [u], [v])\to (x,
[u]),$$ so $D$ is a submanifold of real codimension $4l-2$.

The following lemmas show that why the submanifold $D$ is important.

\begin{lemma} \label{D1} For a nowhere vanishing section $s$ of $E$, if $s$ and $\iota_Es$ are
complex linearly independent, then $([s], [\hat s])\cap
D=\emptyset$. Conversely, if $l=(l_1, l_2)$ is a $\tau-$invariant
section of $P(E)\times_J P({\hat E})$ with $l\cap D=\emptyset$ and
$l_1$ having a lift $s$ to $E$, then $s$ and $\iota_Es$ are complex
linearly independent.
\end{lemma}

\begin{proof} This is clear from definitions.
\end{proof}

\begin{lemma} \label{D2}Any $\tau-$invariant section  of $P(E)\times_J P({\hat E})$ can be deformed to
another one which is transversal to $D$.
\end{lemma}
\begin{proof}
Let $\Gamma$ be a $\tau-$invariant section of $P(E)\times_J P({\hat
E})$. Clearly $\Gamma$ does not intersect $D$ over any fixed point of $J$.
Therefore there is a closed invariant neighborhood ${\mathcal V}$ of
the fixed points set of $J$ over which $\Gamma$ does not intersect
$D$. Let ${\mathcal V}_0$ be the interior of ${\mathcal V}$.
 Away from ${\mathcal V}$, the involution $\tau$ acts freely. Let $P'$
be the quotient of $P(E)\times_JP({\hat E})$ over $J-{\mathcal V}_0$.
Then $P'$ is a $P^{2l-1}\times P^{2l-1}-$bundle over
$J'=(J-{\mathcal V}_0)/\tau$. Since $D$ is $\tau-$invariant,
$D'=P'\cap D/\tau$ is a submanifold of $P'$, in fact a
$P^{2l-1}-$subbundle over $J'$. Since $\Gamma$ is $\tau-$invariant,
it induces  a section $\Gamma'$ over $J'$. Notice that
$\Gamma'|_{\partial J'}$ does not intersect $D'$.  By (the ordinary)
transversality applied to the submanifold $D'\subset P'$ and the map
$\Gamma':J'\to P'$, we can deform $\Gamma'$ to another section
$\Gamma''$ such that $\Gamma''|_{\partial J'}=\Gamma'|_{\partial J'}$
and $\Gamma''$ is transverse to $D'$. The pull back of $\Gamma''$,
together with $\Gamma|_{\mathcal V}$, forms a section of
$P(E)\times_J P({\hat E})$, which is a deformation of $\Gamma$ and
transversal to $D$.

\end{proof}

\begin{theorem} \label{split}  Suppose $J$ has (real) dimension
$k$ and  $E$ is a complex rank $2l$ quaternionic bundle over $J$. If
$4l\geq k+3$ then $E$ splits as $\underline {\mathbb H} \oplus E'$.

\end{theorem}

\begin{proof}
By Lemmas \ref{section} and  \ref{D1}, we just need to construct a
$\tau-$invariant section $l=(l_1, l_2)$ of $P(E)\times_J P({\hat
E})$ with $l\cap D=\emptyset$ and such that $l_1$ has a lift to $E$.


Let $s_0$ be a nowhere vanishing section of $E$. Such a section
exists as $4l\geq k+1$.
  By Lemma \ref{D2}, we can deform
the $\tau-$invariant section $([s_0], [\hat s_0])$ of $P(E)\times_J
P({\hat E})$ to obtain a $\tau-$invariant  section $(l_1, l_2)$
 which is transversal to the $\tau-$invariant
submanifold $D$.

The complex line field $[s_0]$ of $E$ is trivialized by $s_0$. Since
deformations of a trivial complex line field remain trivial, $l_1$
is a trivial complex line field of $E$ as well. In particular, $l_1$
lifts to a nowhere vanishing section $s$ of $E$.

 Since
the dimension of $J$ is  $k$, the dimension of $D$ is equal to
$k+4l-2$, and the dimension of $P$ is equal to $8l-4+k$.
 It follows from the  assumption $4l\geq k+3$ that,
$$\hbox{dim}D+\hbox{dim}l_1=(4l-2+k)+k\leq 4l-2+4l-3+k=8l-5+k=\hbox{dim} P-1.$$
 As $l_1$ is transverse to $D$,  $l_1$ does not intersect $D$.
Therefore $s$ is the required section of $E$.
\end{proof}

\begin{remark} The condition $4l\geq k+3$ in Theorem \ref{split} is sharp since, as mentioned,
there is a non-trivial rank 2 quaternionic bundle over $T^2$.
 On the other hand, it follows from Theorem \ref{split} that  any quaternionic bundle
over $T^1$ is trivial, which is already proved in \cite{FKMM}.
\end{remark}

\begin{corollary}\label{split2}
Suppose $J$ has dimension $4n-\mu$ with $0\leq \mu\leq 3$ and $E$ is
a quaternionic bundle over $J$ with rank $2m\geq 2n$. Then $E$
splits as $Q\oplus \underline{\mathbb H}^{m-n}$ for some rank $2n$
quaternionic bundle $Q$.
\end{corollary}

Next we  give two types of  local trivializations. We first deal
with a sufficiently small invariant disk containing only one fixed
point.

\begin{lemma} \label{near} Any quaternionic bundle is trivial near a fixed point.

\end{lemma}

\begin{proof} Let $U$ be an invariant disk containing only one fixed point. Consider a nowhere vanishing
 section
$s$ over $U$. Since $\iota_E s$ and $s$ are complex linearly
independent over the fixed point, by possibly shrinking $U$ we can
assume they remain so in $U$. Now apply Lemma \ref{section} and
repeat this process.

\end{proof}

Next we treat certain invariant  sets away from fixed points.

 \begin{lemma} \label{away}
 Let $E$ be  complex rank
$2l$ quaternionic bundle  over $J$. Suppose $W$ is a subset of $J$
such that $W$ does not intersect $\iota_J W$ and  $E$ is trivial
over $W$ as a complex vector bundle. Then $E$ is isomorphic to
$\underline {\mathbb H}^l$ over $W\coprod \iota_J W$.
\end{lemma}
\begin{proof}  Let $\alpha:E|_W\to {\mathbb H}^{l}$ be a complex
trivialization of $E$ over $W$. Then
$$\alpha_J:E|_{\iota_JW}\to {\mathbb H}^l, \quad u\to -\alpha (\iota_Eu){\bf j}$$ is
a trivialization of $E$ over $\iota_JW$ and  is complex linear. As
it is assumed that
 $W\cap W'=\emptyset$, $\alpha\coprod \alpha_J$ is a trivialization of $E$
   as a complex bundle over
$W\coprod W'$. Moreover, it is a trivialization of $E$ as a
quaternionic bundle,  since for $u\in E|_{\iota_JW}$, we have,
$$\alpha_J(u){\bf j}=[-\alpha(\iota_E u){\bf j}]{\bf j}=(\alpha\circ \iota_E)(u),$$ and for $v\in E|_W$,
we have
$$(\alpha_J\circ \iota_E)(v)=- \alpha(-v){\bf j}=\alpha(v){\bf j}. $$


\end{proof}

\begin{proposition}\label{cover} For
any quaternionic bundle $E$ over $J$, there is an equivariant
covering  of $J$ such that $E$ is trivial over each open set as a
quaternionic bundle.

\end{proposition}

\begin{proof} For each fixed point $x_i$ of $J$, by Lemma \ref{near} there exists
an open invariant neighborhood $U_i$ of $x_i$ such that $E$ is
trivial over $U_i$ as a quaternionic bundle. Let $V_i$ be a smaller
closed invariant neighborhood of $x_i$ which is  contained in $U_i$.
Let ${\mathcal V}_0$ be the union of the $V_i$. Then  $J-{\mathcal
V}_0$ is invariant and is covered by  disk  pairs $(W_j, \iota_J
W_j)$, where for each $j$,  $W_j$ is a disk contained in
$J-{\mathcal V}_0$ and $W_j\cap \iota_J W_j=\emptyset$. Then, for
each $j$,  $E$ is trivial as a complex bundle over the  disk $W_j$,
and hence trivial as a quaternionic bundle over the invariant open
set $W_j\coprod \iota_J W_j$  by Lemma \ref{away}. Now the $U_i$ and
the $W_j\coprod \iota_J W_j$ form a required covering.
\end{proof}

\begin{example} To illustrate Proposition \ref{cover} we describe an explicit covering of $T^2$.
 Write $T^2$ as $S^1\times S^1$. Cover the p--th $S^1$ by four
disks $A^p_1, A^p_2, B^p_1, B^p_2$. $A^p_1$ and $ A^p_2$ are
disjoint invariant disks around the two fixed points, and they are
called the type $A$ disks; $B^p_1$  and $B^p_2$ are disjoint and
interchanged by the involution, and they are called the type $B$
disks. Consider the union of the
 products of the disks where the type of each factor is fixed.
 Each union consists of $4=2^2$ products of disks.

And there are $4=2^2$ such unions, $${\mathcal U}_{AA}, {\mathcal
U}_{AB}, {\mathcal U}_{BA}, {\mathcal U}_{BB},$$ which form a
covering of $T^2$. Since the involution takes a disk to a disk of
the same type, each union is an invariant subset and so the covering
is equivariant.

We claim that $E$ is trivial as a quaternionic bundle over each
union.
 The $4$ products of $A-$disks in ${\mathcal U}_{AA}$,
called $U_1, ..., U_4$, are disjoint invariant disks in $T^2$, each
containing precisely one fixed point. In particular, by Lemma
\ref{near}, $E$ is trivial over ${\mathcal U}_{AA}$ as a
quaternionic bundle if the $A-$disks are sufficiently small. To show
that $E$ is trivial as a quaternionic bundle over
 each of the remaining 3 unions, we notice that
 each union is  a disjoint union of pairs of product of  disks
interchanged by the involution $\iota_{T^2}$. This is because that
two products are disjoint if and only if some factors are disjoint,
 and distinct disks of the
same type are disjoint. Now apply Lemma \ref{away}.
\end{example}


We could similarly present  an explicit equivariant covering of
$T^k$, which might be used to give another calculation of $KQ(T^k)$
in \cite{FK}, and to extend the classification in \cite{FKMM} to all
$T^k$.

\begin{theorem} \label{monomorphism} Suppose $E$ is a rank $2l$ quaternionic bundle over a compact $J$. Then
there is a  $Pin(2)-$equivariant monomorphism from $E$ to
$\underline{\mathbb H}^{m}$ for $m=l+[\frac {k+2}{ 4}]$. Here $[x]$
denotes the largest integer bounded by $x$ from above.
\end{theorem}
\begin{proof} Consider an equivariant covering  $\{U_i, W_j\coprod \iota_J W_j\}_{i, j}$ as in Proposition \ref{cover}.
We can assume this covering is finite as $J$ is compact. By possibly
shrinking the $U_i$ we can assume that $U_i\cap U_j=\emptyset$ for
$i\ne j$. Therefore we can trivialize $E$ as a quaternionic bundle
over ${\mathcal  U}=\coprod_i U_i$.  Since $m=l+[{k+2\over 4}]\geq
l$, we  can view this trivialization as a $Pin(2)-$equivariant
monomorphism $\Phi_0$ from $E$ to $\underline {\mathbb H}^{m}\supset
\underline {\mathbb H}^{l} $ over ${\mathcal U}$.

Let ${\mathcal W}_0={\mathcal U}$, and for $j\geq 1$, let
$${\mathcal
W}_j={\mathcal W}_{j-1}\cup (W_j\coprod \iota_J W_j).$$ We will
argue by induction on $j$. Suppose the $Pin(2)-$equivariant
monomorphism $\Phi_j$ has been defined over ${\mathcal W}_{j}$. Over
$W_{j+1}$, fix a complex trivialization
$$\Psi_{j+1}:E|_{W_{j+1}}\to {\mathbb H}^l.$$

Let $K={\mathcal  W}_j\cap W_{j+1}$. Then for each $x\in K$,
$$\phi_{j+1}=\Phi_j\circ \Psi_{j+1}^{-1}:{\mathbb H}^l\to E|_x\to {\mathbb H}^m$$
is a complex monomorphism, and hence a point in the complex Stieffel
manifold $V_{2m, 2l}$ of linear maps from ${\mathbb C}^{2l}$ to
${\mathbb C}^{2m}$ of rank $2l$.

The space $V_{2m, 2l}$ naturally lies inside ${\mathbb C}^{2l\times
2m}$. We can use a partition of unity to extend $\phi_{j+1}$ as a
map from $K$ to ${\mathbb C}^{2l\times 2m}$ to a map
$$\tilde \phi_{j+1}: W_{j+1}\to {\mathbb
C}^{2l\times 2m}.$$
  We would like the extension to actually lie in $V_{2m, 2l}$.
This is achieved by a transversality argument.
 The  complement of $V_{2m, 2l}$ is stratified by linear maps of lower ranks.
The  stratum $S_{2l-b}$  with rank $2l-b$ is a submanifold with real
codimension $2(2m-2l+b)b$. We assume the extension $\tilde
\phi_{j+1}$ is transversal to all $S_{2l-b}, 1\leq b\leq 2l$.

The stratum $S_{2l-1}$ with rank $2l-1$ has the smallest
codimension, which is
\[
2(2m-2l+1)=4[{k+2\over 4}]+2= \left\{\begin{array}{ll} k+2, & \hbox{
if } k\equiv 0
\pmod 4\\
        k+1, & \hbox{ if } k\equiv 1 \pmod 4\\
        k+4, & \hbox{ if } k\equiv 2 \pmod 4\\
        k+3, & \hbox{ if } k\equiv 3 \pmod 4,\\
        \end{array}
\right. \] because  $m=l+[{k+2\over 4}]$. Since the dimension of
$W_{j+1}$ is $k$,  $\tilde \phi_{j+1}$  misses each $S_{2l-b}, 1\leq
b\leq 2l$.

Now at each point $y\in W_{j+1}$,  $$\tilde \phi_{j+1}\circ
\Psi_{j+1}: E|_y \to {\mathbb H}^l \to {\mathbb H}^m$$ is a complex
monomorphism, and it agrees with $\Phi_{j}$ over $K$. As in Lemma
\ref{away} we can canonically extend it  $Pin(2)-$equivariantly to
$\iota_J W_{j+1}$. Since $\Phi_j$ is assumed to be
$Pin(2)-$equivariant, the extension also agrees with $\Phi_j$ over
${\mathcal W}_j\cap \iota_J W_{j+1}$. Thus we obtain a
$Pin(2)-$equivariant monomorphism
$$\Phi_{j+1}:E|_{{\mathcal W}_{j+1}}\to \underline{\mathcal
H}^m.$$
\end{proof}

\begin{example}
According to Theorem \ref{monomorphism} any rank $2$ quaternionic
bundle over $T^4$ can be embedded into $\underline{\mathbb H}^2$,
since $1+[\frac{4+2}{4}]=2$. This can be also proved using
\cite{FKMM}. Indeed, it is shown there that
  any rank $2$ bundle over $T^4$ is of the form
$E=\underline{\mathbb H}(S)$, where $S$ is a signed invariant finite
set of $T^4$ and $\underline{\mathbb H}(S)$ is obtained from
$\underline{\mathbb H}$ by a canonical spinor twisting around $S$.
Moreover, if
 $ E'=\underline{\mathbb H}(-S)$, then  $E\oplus E'=\underline{\mathbb H}^2$.
\end{example}

\section{ Stable homotopy and unoriented bordism Seiberg-Witten
invariants} \label{homotopy and bordism}

In this section  $M$ is a closed oriented smooth 4-manifold
and $c$ is a Spin$^c$ structure.
We first review the  stable homotopy Seiberg-Witten invariants.
Then we construct the unoriented bordism Seiberg-Witten invariants.

\subsection{Seiberg-Witten equations}

Let $S^0$ and $ S^1$
be the  spinor bundles associated to $c$.
 The determinant line bundles  $\hbox{det}_{\mathbb C} S^0$ and
$\hbox{det}_{\mathbb C} S^1$ are isomorphic. Denote  the Hermitian
line bundle by $L_c$. Fix a  Hermitian connection $A_0$ on $L_c$.
Let ${\mathcal H}^1(M,{\mathbb R})$ be the space of harmonic
1-forms, and
 consider the affine space ${\mathcal  A}_0$ of Hermitian connections on $L$ of
  the form $A=A_0+a{\bf i}$ for
$a\in {\mathcal  H}^1(M,{\mathbb R})$

Let ${\mathcal  H}^0(M,U(1))$ be the group of harmonic maps from $M$
to $U(1)$. Fix a base point $x_0\in M$ and let ${\mathcal
H}_0^0(M,U(1))$ be the subgroup consisting of the harmonic maps
sending $x_0$ to the identity. Then ${\mathcal H}^0(M,U(1))$ is the
product of $U(1)$ and ${\mathcal H}_0^0(M, U(1))$, where $U(1)$ is
the subgroup of constant maps.

 Consider the gauge group action of
$g\in {\mathcal  H}^0(M,U(1))$ on $A\in {\mathcal  A}_0$ and let $J$
be the quotient of ${\mathcal  A}_0$ by ${\mathcal  H}_0^0(M,
U(1))$. Then $J$ is identified with  the
 quotient of $H^1(M;{\mathbb R})$ by
$H^1(M;2{\mathbb Z})$, and thus a torus of dimension $b_1$.

Let $p_M:{\mathcal A}_0\times M\to M$ be the projection map.
${\mathcal H}^{0}(M, U(1))$ acts on  the bundle $p_M^* (S^0\oplus
S^1)$ on ${\mathcal A}_0\times M$ by
$$
g:(A, x, q)\longrightarrow (g(A), x, qg).
$$
The action of the subgroup $U(1)$ coincides with the action coming
from the complex structure.

Consider the two infinite dimensional complex vector bundles over
$J$:
$$
\begin{array}{ll}\tilde V^0&=({\mathcal  A}_0\times \Gamma(p_M^*S^0))/{\mathcal
H}_0^0(M, U(1)),\cr \tilde V^1& =({\mathcal A}_0\times
\Gamma(p_M^*S^1))/{\mathcal  H}_0^0(M, U(1)),\cr
\end{array}
$$
and the smooth family of $U(1)-$equivariant Dirac operators
$\{D_a\}_J$. According to  the Atiyah-Singer index theorem  we have
the formula
\begin{equation}\label{index}
\hbox{rank}_{\mathbb C}\hbox{Ind}\{D_a\}_J=-\sigma(M)/8.
\end{equation}

 Consider as well the  two trivial infinite
dimensional real vector bundles over $J$:
$$\begin{array}{ll}\tilde W^0& = J\times d^{*}(\Gamma(\Lambda^2TM^*))\subset J\times
\Gamma(TM^*),\cr \tilde W^1&=J\times \Gamma(\Lambda^{+}TM^*),\cr
\end{array}$$ with trivial $U(1)-$action and the operator $d^+$,
which is the self-dual part of $d$.

With this set up, the Seiberg-Witten equations are then  a
$U(1)-$equivariant bundle map $\tilde f_{SW}$ between the infinite
dimensional bundles $\tilde V^0\oplus \tilde W^0$ and $\tilde
V^1\oplus \tilde W^1$,
 which, at a point $a\in J$, is of the form
$$\tilde f_{SW}|_a(s, b)=(D_{a}s+\frac{1}{2}C(b)s{\bf i}, d^+b+s{\bf i}\bar s).$$
Here $C:TM^*\otimes S^0\to S^1$ is the Clifford multiplication, and
$s{\bf i}\bar s$ is a natural algebraic map from $S^0$ to
$\Lambda^{+}TM^*$ (see e.g. \cite{L}).

\begin{remark}\label{independent}When restricted to $0\oplus \tilde W_0$, $\tilde f_{SW}$ is the
linear (embedding)
 sending $(0,b)$ to $(0, d^+b)$ at each $a\in J$, in particular,
it is independent of $a\in J$.
\end{remark}

Now let us suppose $c$ is a reducible Spin$^c$ structure. Then $L_c$
is a trivial bundle.  Coming with a reduction of $c$ to a spin
structure are
 the involution on $J$ and the enlarged $Pin(2)$ symmetry
 of $\tilde f_{SW}$, which we explain briefly now.

Fixing a trivial connection $A_0$ on $L_c$ corresponding to the spin
reduction,
 and let $D_0$ be the associated
Dirac operator. Consider the involution $\iota_{{\mathcal  A}_0}$ on
${\mathcal  A}_0$ by sending $a$ to $-a$. We lift $\iota_{{\mathcal
A}_0}$ to the bundle $p_M^* (S^0\oplus S^1)$ by
$$
\iota_S:(a, x, s)\longrightarrow (-a, x, s{\bf j})$$ to make both
$p_M^*S^0$ and $p_M^*S^1$ quaternionic bundles. When dividing out by
${\mathcal H}_0^0(M, U(1))$, $\iota_{{\mathcal A}_0}$ induces the
standard involution $\iota_J$ on  the torus $J$. Furthermore,
$\iota_S$ induces the anti-complex lifts $\iota_{\tilde V^0}$ and
$\iota_{\tilde V^1}$ of $\iota_J$ on the bundles ${\tilde V}^0$ and
${\tilde V}^1$, which make them quaternionic bundles. Therefore
$Pin(2)$ acts on both $\iota_{\tilde V^0}\oplus \iota_{\tilde W^0}$
and $\iota_{\tilde V^1}\oplus \iota_{\tilde W^1}$ by Remark
\ref{real}.

 It is well-known that  $D_0$ is a ${\mathbb H}-$linear
operator. From which it is not hard to see that  the family of
linear operators $\{D_a\}_J$ is $Pin(2)-$equivariant. Hence the
index bundle of $\{D_a\}_J$, $\hbox{Ind}\{D_a\}_J$, lies in $KQ(J)$.
For the $Pin(2)-$equivariance of the remaining terms of the  map
$\tilde f_{SW}$ we refer to \cite{L}.

An important property of the Seiberg-Witten equations is that
$\tilde f_{SW}$ is proper, which  implies that we formally have a
$Pin(2)-$equivariant map between the two infinite dimensional sphere
bundles and thus an element $\tilde f_{SW}\in [S(\tilde V^0\oplus
\tilde W^0), S(\tilde V^1\oplus \tilde W^1)]^{Pin(2)}$. Here we
follow the notations in \cite{FKMM}: For $G=U(1)$ or $Pin(2)$,
$\hbox{Map} (S(V), S(V'))^G$ denotes the set of $G-$equivariant maps
between the sphere bundles of the $G-$equivariant bundles $V$ and
$V'$. Its quotient divided by the $G-$equivariant homotopy is
denoted by $[S(V), S(V')]^{G}$. To understand this element more
explicitly, we need the technique of finite dimensional
approximations (initiated in \cite{F1}), which leads to the
construction of the stable cohomotopy Seiberg-Witten invariants of
Bauer and Furuta.

\subsection{Finite dimensional approximations}

First of all we need the notion of an admissible
$Pin(2)-$equivariant triple.

\begin{definition} \label{triple}A
$Pin(2)-$equivariant triple is a triple
$$(S( V_0\oplus W_0),
S(V_1\oplus  W_1), f)$$ consisting of

\begin{enumerate}
\item finite dimensional quaternionic vector bundles $V_0$ and $V_1$
over $J$,

\item finite dimensional trivialized real vector bundles $W_0$ and $ W_1$
over $J$,

\item a $Pin(2)-$equivariant map
$$f:S(V_0\oplus  W_0) \longrightarrow S(V_1\oplus  W_1).$$
\end{enumerate}
A $Pin(2)-$equivariant triple is called admissible if,  when
restricted to $0\oplus W_0$, $f$
 is independent of $a\in J$.
In the same way we define
admissible $U(1)-$equivariant triples.

\end{definition}

Notice that $f$ maps $0\oplus W_0$ to $0\oplus W_1$ due to
equivariance. And since $W_0$ and $W_1$ are trivialized, it makes
sense to require the restriction of $f$ to $0\oplus W_0$ be
independent of $a\in J$.

 Let us recall the stabilization process.
Given a $Pin(2)-$equivariant  triple
$(S( V_0\oplus  W_0), S(V_1\oplus  W_1),f)$,
 a quaternionic vector bundle $V$ and
a trivial real vector bundle $W$, define
 \begin{equation}
\label{stabilization}\hat V_0=V_0\oplus V,\quad \hat V_1=V_1\oplus
V,\quad \hat W_0=W_0\oplus W, \quad \hat W_1=W_1\oplus W.
\end{equation}

Recall that the sphere bundles $S(V_i\oplus W_i)$ and $S(V\oplus W)$
can be joined fibrewisely to form the sphere bundle of the direct
sum $\hat V_i\oplus \hat W_i$ by the formula $(1-t) a_i+t a$ for
$t\in [0,1]$, and similarly two $Pin(2)-$equivariant maps
$$f:S(V_0\oplus
W_0)\longrightarrow
S(V_1\oplus W_1) \quad  \hbox{  and   }\quad  g:S(V\oplus W)\longrightarrow
S(V\oplus W).$$
can be joined to a $Pin(2)-$equivariant map
$$j(f,g):S( \hat V_0\oplus \hat W_0) \longrightarrow  S(\hat V_1\oplus \hat W_1).$$
Thus by taking  the join  with the identity on $S(V\oplus W)$, we have the
stabilization map between the $Pin(2)-$equivariant triples
$$ (S( V_0\oplus  W_0), S(V_1\oplus  W_1), f)\longrightarrow
(S( \hat V_0\oplus \hat W_0), S(\hat V_1\oplus \hat W_1),
j(f, id)).$$

It is easy to check that the join of two homotopies is a homotopy.
Two triples are called stable homotopic if they become homotopic
under  stabilization. This is an equivalence relation. We call an
equivalence class of triples a $Pin(2)-$equivariant stable homotopy
class. Clearly the join of two admissible  triples is still
admissible.


A finite dimensional approximation to the Seiberg-Witten equations
associated to a spin reduction of $c$
 is a $Pin(2)-$equivariant triple $(S( V_0\oplus  W_0), S(V_1\oplus  W_1),
f_{SW})$ such that
\begin{equation} \label{approximation}
\begin{array}{ll}
[V_0]-[V_1]&=\hbox{ind}\{D_a\}_J\in KQ(J),\cr
[W_0]-[W_1]&=[-\underline{\mathbb R}^{b^+}]\in KO(J).\cr
\end{array}
\end{equation}
And $f_{SW}$ is an approximation of $\tilde f_{SW}$ in an
appropriate sense, which we do not specify as it will be irrelevant
for us (see Proposition \ref{constant}). An admissible finite
dimensional approximation is one such that $f_{SW}$ is independent
of $a\in J$ when restricted to $S(0\oplus W_0)$.

There are many such finite dimensional approximations, all of which
are related via the stabilization process. More precisely, it was
shown in \cite{BF} and \cite{F} that any two finite dimensional
approximations become homotopic under  stabilization, and moreover, the
 homootpy itself is well-defined up to homotopy.
 Notice that it
is pointed out in Remark \ref{independent} that the SW map $\tilde
f_{SW}$ is admissible in the sense that it is independent of $a\in
J$ when restricted to $0\oplus \tilde W^0$.
 Indeed it is further
shown in \cite{B} and  \cite{FL} that it can be assumed that the
finite dimensional approximations are admissible.

Therefore  there is  a well-defined $Pin(2)-$equivariant stable
homotopy class of admissible triples.
 Furthermore, this $Pin(2)-$equivariant stable homotopy class
only depends on the oriented diffeomorphism type of $M$ and the
reducible Spin$^c$ structure $c$ and the spin reduction $\nu$  (see
\cite{F}). Thus, we can write this $Pin(2)-$equivariant stable
homotopy class as
 $\mathcal  {SW}(M,c,\nu)$ and call it the $Pin(2)-$stable cohomotopy Seiberg-Witten
invariant.

\subsection{Unoriented bordism SW invariants}
In this subsection we  construct an unoriented bordism SW invariant
of a reducible Spin$^c$ structure  when $b^+\geq 2$ following
\cite{FKMM}. Our invariants are simpler, living in the unoriented
bordism group rather than the richer Pin bordism group as in
\cite{FKMM}. But this would be sufficient for our purpose.

On the other hand, we only need the assumption $b^+\geq 2$ rather
than $b^+\geq b_1+2$. Being able to weaken the assumption on $b^+$
is crucial for us. This is achieved by adding the admissibility as
in Definition \ref{triple}.


The construction of  the unoriented bordism SW invariant is given in
several steps.

\subsubsection{The construction of $\gamma'$}
Let $$\hbox{Map}_{adm}(S( V_0\oplus  W_0), S(V_1\oplus
W_1))^{U(1)}$$ be the space of $U(1)-$equivariant maps between the
sphere bundles which are admissible. Finite dimensional
approximations to the Seiberg-Witten equations give rise to such
objects.

Given $f_0$ and $f_1$ in $\hbox{Map}_{adm}(S( V_0\oplus  W_0),
 S(V_1\oplus  W_1))^{U(1)}$, we
can view them as maps from $S(V_0\oplus W_0)$ to $V_1\oplus W_1$.
Let $F(f_0, f_1)$ be the set of $U(1)-$equivariant paths
$$\tilde f_t:S(V_0\oplus W_0)\times [0,1]\longrightarrow V_1\oplus W_1$$
 connecting $f_0$ and $f_1$, and  satisfying
\begin{enumerate}
\item the restriction of $\tilde f_t$ to $S(0\oplus W_0)\times [0,1]$,
 which is mapped to $S(0\oplus W_1)$ due to equivariance,
 does not vanish and  is independent of $a\in J$, i.e.
 there is a map  $\xi$  from $S({\mathbb R}^{\hbox{rank}_{\mathbb R}W_0})\times [0,1]$ to ${\mathbb
R}^{\hbox{rank}_{\mathbb R}W_1}-0$ such that $\tilde f_t=\xi$ on
$S(0\oplus W_0)|_a\times [0,1]$ for any $a\in J$.

\item $\tilde f_t$ is transverse to the zero section.
\end{enumerate}

From now on we assume that in this section
\begin{equation}\label{W}
\hbox{rank}_{\mathbb R}W_1-\hbox{rank}_{\mathbb R}W_0\geq 2.
\end{equation}
This corresponds to $b^+\geq 2$.

\begin{lemma} \label{nonempty} $F(f_0, f_1)$ is non-empty.
\end{lemma}
\begin{proof}
 The existence of a $\tilde f_t$ is shown by three steps.

{\bf Step 1}. Since the fibers of $V_1\oplus W_1$ are linear spaces,
we can use
simply a linear combination  to construct a $U(1)-$equivariant map
$\Omega$ from $S(V_0\oplus W_0)\times [0,1]$ to $V_1\oplus W_1$
connecting $f_0$ and $f_1$.

{\bf Step 2}. Both $f_0$ and $f_1$  are assumed to be independent of
$a\in J$ when restricted to $S(0\oplus W_0)$. By the assumption
(\ref{W}), every two maps from $S({\mathbb R}^{\hbox{rank}_{\mathbb
R}W_0})$ to ${\mathbb R}^{\hbox{rank}_{\mathbb R}W_1}-0$ are
homotopic. In particular, the restrictions of $f_i$, as maps from
$S(0\oplus W_0)$ to $0\oplus(W_1-0)$ are
 homotopic to each other
  through a homotopy which is independent of $a\in J$.
 Therefore we can perturb
$\Omega$ near $S(0\oplus W_0)\times [0,1]$ but away from
$S(V_0\oplus W_0)\times (0\coprod 1)$, by an equivariant partition
of unity, to an $U(1)-$equivariant homotopy $\Omega'$ connecting
$f_0$ and $f_1$, and such that $\Omega'$ satisfies (1).


{\bf Step 3}. Now the zero set of $\Omega'$ is away from the closed
subsets
$$S(0\oplus W_0)\times [0,1] \quad \hbox{and} \quad S(V_0\oplus W_0)\times
(0\coprod 1).$$ Thus it has a neighborhood with the same property.
In particular  $U(1)$ acts freely on such a neighborhood. Hence we
can further perturb $\Omega'$ equivariantly in such a neighborhood
to make it transverse to the zero section. The new perturbation is
then a homotopy connecting $f_0$ and $f_1$, and satisfies both (1)
and (2).
\end{proof}

Given  $\tilde f\in F(f_0, f_1)$, denote
 the zero set of $\tilde f^{-1}(0)$ by $\hat {\mathcal  M}$.
Then $\hat {\mathcal  M}$ is a smooth,  closed submanifold of
$S(V_0\oplus W_0)\times [0,1]$. Let $\hat B$ be the complement of
$S(0\oplus W_0)$ in $S( V_0\oplus W_0)$. Then  $\hat {\mathcal M}$
actually lies in $\hat B\times (0,1)$, so it is itself a closed
smooth manifold. Denote the $U(1)-$quotient $\hat {\mathcal M}/U(1)$
by ${\mathcal M}$. Since $U(1)$ acts freely on $\hat B$, the
quotient
 ${\mathcal M}$
 is also a closed smooth manifold.
 Moreover, we can view $\tilde f_t$ as a section of
 the bundle
 $${\mathcal E}=(S(V_0\oplus W_0)\times [0,1])\times_J (V_1\oplus
 W_1)$$ over $S(V_0\oplus W_0)\times [0,1]$. In particular,
 the dimension of ${\mathcal M}$ is easily seen
to be
\begin{equation}\label{dimension}
\dim J-1+2\hbox{rank}_{\mathbb C}V_0-2\hbox{rank}_{\mathbb C}V_1+
\hbox{rank}_{\mathbb R}W_0-\hbox{rank}_{\mathbb R}W_1.
\end{equation}

\begin{lemma}\label{bordism}
The unoriented bordism class of ${\mathcal  M}$ does not depend on
the choice of $\tilde f\in F(f_0, f_1)$.
\end{lemma}

\begin{proof}
Given  $(\tilde f_{t})_0$ and $(\tilde f_{t})_1$ in $F(f_0, f_1)$,
we can construct a homotopy
$$\tilde f_{t, s}:(S(V_0\oplus W_0)\times [0,1])\times [0,1]\longrightarrow V_1\oplus W_1$$
 such that
\begin{enumerate}

\item  $\tilde f_{t, 0}=(\tilde f_{t})_0$ and $\tilde f_{t,
1}=(\tilde f_{t})_1$,

\item the restriction of $\tilde f_{t, s}$ to $S(V_0\oplus W_0)\times
(0\coprod 1)\times s$ is independent of  $s\in [0,1]$, i.e.
\begin{equation} \label{0and1}\begin{array}{ll}
\tilde f_{0, s}=&\tilde f_{0,0}=\tilde f_{0,1}=f_0, \\
\tilde f_{1, s}=&\tilde f_{1,0}=\tilde f_{1,1}=f_1, \\
\end{array}
\end{equation}
 for any $s\in [0,1]$, in  particular, $\tilde f_{t, s}$ does not vanish
 on $$S(V_0\oplus W_0)\times (0\coprod 1)\times [0,1],$$

\item the restriction of $\tilde f_{t, s}$ to $S(0\oplus W_0)\times
[0,1]\times [0,1]$ does not vanish,

\item $\tilde f_t$ is transverse to the zero section.
\end{enumerate}

The zero set of $\tilde f_{t,s}$ then is a compact manifold whose
only boundaries are  $\hat {\mathcal M}_0$ and $\hat {\mathcal
M}_1$. In addition,  $U(1)$ acts freely on it. The smooth
$U(1)-$quotient then provides the desired bordism.

The existence of $\tilde f_{t,s}$ is established in the same way as
that of $\tilde f_t$. We first construct a $U(1)-$equivariant map
$\tilde \Omega$ from $S(V_0\oplus W_0)\times [0,1])\times [0,1]$ to
$V_1\oplus W_1$ connecting   $
 (\tilde f_{t})_0$ and $
 (\tilde f_{t})_1$, and  such that
(2),  or equivalently, (\ref{0and1}), is satisfied.  For example we
could use a linear homotopy.

There are maps
$$\xi_i:S({\mathbb R}^{\hbox{rank}_{\mathbb
R}W_0})\times [0,1]\longrightarrow  {\mathbb
R}^{\hbox{rank}_{\mathbb R}W_1}-0$$ such that $\tilde f_i=\xi_i$ on
$S(0\oplus W_0)|_a\times [0,1]$ for any $a\in J$.
 From the assumption
(\ref{W}), the maps $\xi_i$  are homotopic relative to $S({\mathbb
R}^{\hbox{rank}_{\mathbb R}W_0})\times (0\coprod 1)$. Therefore the
restrictions of $\tilde f_i$ to $S(0\oplus W_0)\times [0,1]$ are
homotopic as maps to  $0\oplus(W_1-0)$
  through a homotopy which is constant on
$S(0\oplus W_0)\times (0\coprod 1)$ (and independent of $a\in J$).
Thus we can perturb $\tilde \Omega$ near
$$S(0\oplus W_0)\times [0,1]\times [0,1]$$
but away from
$$S(V_0\oplus W_0)\times [0,1]\times (0\coprod 1)\quad \hbox{and}
\quad S(V_0\oplus W_0)\times (0\coprod 1)\times [0,1]  $$
 to an $U(1)-$equivariant homotopy $\tilde \Omega'$
 connecting $
 (\tilde f_{t})_0$ and $
 (\tilde f_{t})_1$, and
such that both (2) and (3) are satisfied.

Now the zero set of $\tilde \Omega'$ is away from the closed subsets
$S(0\oplus W_0)\times [0,1]\times [0,1]$ and $S(V_0\oplus W_0)\times
(0\coprod 1)\times [0,1]$. Hence it has a neighborhood $U$ with the
same property. In particular, $U(1)$ acts freely on $U$.   Hence we
can further perturb $\tilde \Omega'$ equivariantly inside $U$ to make it
transverse to the zero section. Notice that $\tilde \Omega'$ is
already transverse to the zero section along $S(V_0\oplus W_0)\times
[0,1]\times (0\coprod 1)$, so the perturbation can be chosen to be
also away from $S(V_0\oplus W_0)\times [0,1]\times (0\coprod 1)$.
The new perturbation is then a homotopy connecting $
 (\tilde f_{t})_0$ and $
 (\tilde f_{t})_1$, and
such that (2), (3) and (4)  are all satisfied.
\end{proof}

Therefore we can make the following definition.

\begin{definition}  \label{gamma'}For  $f_0, f_1$ in $\hbox{Map}_{adm}(S( V_0\oplus  W_0),
S(V_1\oplus  W_1))^{U(1)}$, let $\gamma_{V_0\oplus W_0, V_1\oplus
W_1}' (f_0, f_1)=[{\mathcal M}]^{uo}\in \Omega_n^{uo},$ where $n$ is
given by (\ref{dimension}).

\end{definition}

\subsubsection{Properties of $\gamma'$}
 We now establish a few properties of $\gamma'$.

 Since (\ref{W}) is  invariant under
 stabilization, $\gamma'_{-,-}$ is defined on any stabilization of
the pair of maps $f_0$ and $f_1$.   Furthermore, the bordism class does not
change, as
the join of $\tilde f$ and $id$ has the same zero set as that
of $\tilde f$.

In addition, $\gamma_{V_0\oplus W_0, V_1\oplus W_1}'$ satisfies an
important additivity property.
 Given $\tilde f\in F(f_0, f_1)$ and $\tilde g\in F(f_1, f_2)$, they
 naturally combine  to
an element $\tilde h\in F(f_0, f_2)$, defined by
\[
\tilde h(t)=\left \{\begin{array}{ll}\tilde f(2t), &\hbox{if $0\leq
t\leq {1\over 2}$},\\ \tilde g(2t-1), &\hbox{if $ {1\over 2}\leq
t\leq 1$}.\\
\end{array}
\right.\]
 Clearly the zero set of $\tilde h$ is the disjoint union of
those of $\tilde f$ and $\tilde g$. Therefore
 $\gamma_{V_0\oplus W_0, V_1\oplus W_1}'$ is additive in the following sense:
 \begin{equation}\label{additivity}\gamma_{V_0\oplus W_0, V_1\oplus W_1}'(f_0, f_1)+
\gamma_{V_0\oplus W_0, V_1\oplus W_1}'(f_1, f_2)=\gamma_{V_0\oplus W_0, V_1\oplus W_1}'(f_0, f_2).
\end{equation}

This additivity immediately implies that $\gamma_{V_0\oplus W_0,
V_1\oplus W_1}'$  only depends on the homotopy classes of $f_0$ and
$f_1$. Thus, we can and will from now on regard $\gamma_{V_0\oplus
W_0, V_1\oplus W_1}'$ as a map from $$[S(V_0\oplus W_0), S(V_1\oplus
W_1)]_{adm}^{U(1)}\times [S(V_0\oplus W_0), S(V_1)\oplus
W_1)]_{adm}^{U(1)}$$
 to $\Omega_n^{uo}$.
Obviously the additivity still holds with this new meaning of $\gamma_{V_0\oplus W_0, V_1\oplus W_1}'$

Another  important property of $\gamma'$ is the following.
\begin{proposition} \label{constant} For $Pin(2)-$equivariant sections,  $\gamma'$ is  independent of homotopy classes
of $Pin(2)-$equivariant maps, i.e. the composition
$$\gamma': [S(V_0\oplus W_0), S(V_1\oplus W_1)]_{adm}^{Pin(2)}\times
[S(V_0\oplus W_0), S(V_1)\oplus W_1)]_{adm}^{Pin(2)} \longrightarrow
\Omega_n^{uo}$$ is a constant map.

\end{proposition}
\begin{proof} Consider two $Pin(2)-$equivariant maps $f_0$ and $f_1$ between
the pairs. Notice that, as $U(1)$,  $Pin(2)$ acts freely away from
the $U(1)-$fixed point set
$${\mathcal F}=S(0\oplus W_0)\coprod S(0\oplus W_1).$$ Notice also
that $\iota$ acts freely on the set ${\mathcal F}$ as an involution.
Applying the dimension assumption (\ref{W}) to the quotient
manifolds of ${\mathcal F}/\iota$, we can actually
 construct a $\tilde f_t\in F(f_0, f_1)$ which is $Pin(2)-$equivariant.

Thus $\iota$ is a free involution on ${\mathcal M}$. Let
$p:{\mathcal M}\to {\mathcal M}/\iota$ be the double covering and
$\zeta$ the real line bundle associated to $p$. Then ${\mathcal M}$
is diffeomorphic to the sphere bundle of $\zeta$, hence it bounds
the disk bundle of $\zeta$.
 Therefore
the unoriented bordism class of ${\mathcal M}$ is zero, that is,
 $\gamma_{V_0\oplus W_0, V_1\oplus W_1}'([f_0], [f_1])=0$.

Together with the additivity of $\gamma_{V_0\oplus W_0, V_1\oplus
W_1}'$, we have the proposition.

\end{proof}

\subsubsection{The invariant $e_1(V_0\oplus W_0, V_1\oplus W_1)$}
We first construct a variation of $\gamma'$, $\gamma_{V_0\oplus W_0,
V_1\oplus W_1}$, whose input is a single element, rather than a
pair,  in $[S(V_0+W_0), S(V_1\oplus W_1)]_{adm}^{U(1)}$.

  Consider constant maps in $\hbox{Map} (S(V_0+W_0),
S(V_1\oplus W_1))^{U(1)}$. By the $U(1)-$equivariance, they must
land  in $S(0\oplus W_1)$.
By the assumption (\ref{W}) all such maps are homotopic. Let $[f_0]$
be this unique homotopy class of constant maps. For any $[f]\in
[S(V_0+W_0), S(V_1\oplus W_1)]_{adm}^{U(1)}$ we define
$$\gamma_{V_0\oplus W_0, V_1\oplus W_1} ([f])=\gamma_{V_0\oplus W_0, V_1\oplus W_1}'([f_0], [f]).$$
$\gamma_{V_0\oplus W_0, V_1\oplus W_1}$ is also invariant under
stabilization since the join of a constant map $f_0$ and $id$ is
itself homotopic to a constant map (just observe that the join of a
point and a sphere is a disk). By Proposition \ref{constant} and
(\ref{additivity}), $\gamma$ takes a constant value on
$Pin(2)-$equivariant sections.

\begin{definition}\label{e_1} Suppose $[S(V_0\oplus W_0), S(V_1\oplus
W_1)]_{adm}^{Pin(2)}\ne \emptyset$.  We write the constant image of
$\gamma$ on $[S(V_0\oplus W_0), S(V_1\oplus W_1)]_{adm}^{Pin(2)}$ in
$\Omega_n^{uo}$ as
$$e_1(V_0\oplus W_0, V_1\oplus W_1).$$
\end{definition}

Since $\gamma_{V_0\oplus W_0, V_1\oplus W_1}$ is invariant under
stabilization, $e_1$ satisfies the stabilization property:
\begin{equation} \label{stability} e_1(V_0\oplus W_0, V_1\oplus W_1)=e_1(
\hat V_0\oplus \hat W_0, \hat V_1\oplus \hat W_1),
 \end{equation}
where $\hat V_i$ and $\hat W_i$ are given as in
(\ref{stabilization}).

\subsubsection{The unoriented bordism SW invariant}  Now let $M$ be a spin manifold
with $2\chi+3\sigma=0$ and $b^+\geq 2$. Let $c$ be a reducible
Spin$^c$ structure together with a spin reduction $\nu$. Then a
finite dimensional approximation $(S(V_0\oplus W_0), (V_1\oplus
W_1), f_{SW})$  is $Pin(2)-$equivariant and can be chosen to be
admissible. In particular,
 $[S(V_0\oplus W_0), S(V_1\oplus W_1)]_{adm}^{Pin(2)}$ is nonempty. Thus we can make the
following definition, in view of (\ref{stability}).

\begin{definition}  Suppose  $M$ is a manifold with $b^+\geq 2$ and
 $c$ is a reducible Spin$^c$ structure on $M$ together with a spin reduction $\nu$.
The unoriented bordism Seiberg-Witten invariant $SW^{uo}(M, c)
:\mathcal  {SW}(M,c,\nu)\longrightarrow \Omega_{n(M,c)}^{uo}$ is
defined to be
 $$
SW^{uo}(M, c)=e_1(V_0\oplus W_0, V_1\oplus W_1)
$$
where $V_0, W_0, V_1,W_1$ arise from an admissible finite
dimensional approximation of the Seiberg-Witten equations associated
to $(c,\nu)$ and $n(M,c)$ is given by (\ref{dimension}).

\end{definition}
It turns out  $SW^{uo}(M, c)$  is independent of $\nu$, and is
 an invariant of
the oriented diffeomorphism type of $M$ and $c$. This is because for
different spin reductions the admissible finite dimensional
approximations are still $U(1)-$equivariantly stably homotopic.
Hence they will have the same $\gamma$ invariant due to the
stability of $\gamma$.

Notice that in this case $J=T^{b_1}$ and $V_0, W_0, V_1,W_1$
 satisfy (\ref{approximation}).
Therefore, we have
\begin{equation} \label{dimension2}
n(M,c)=b_1-1-\frac{\sigma}{4}-b^+=\frac{4b_1-4-5b^++b^-}{4},
\end{equation}
by (\ref{dimension}), (\ref{approximation}, (\ref{index}), and
\begin{equation}\label{signature}
\sigma=b^+-b^-.
\end{equation}
 Since
\begin{equation} \label{euler}
\chi=2-2b_1+b^++b^-,
\end{equation}
we have
\begin{equation}\label{2euler+3signature}
2\chi+3\sigma=4-4b_1+5b^+ -b^-.
\end{equation}
Recall that the SW moduli space of the reducible Spin$^c$ structure
$c$  is
\begin{equation}\label{swdimension} -\frac{2\chi+3\sigma}{4}+\frac{c_1(L_c
)^2}{4}=-\frac{2\chi+3\sigma}{4}, \end{equation}
 as $L_c$ is a
trivial bundle.

Comparing (\ref{dimension2}), (\ref{2euler+3signature}) and
(\ref{swdimension}), we find that $n_{M,c}$ agrees with the
dimension of the SW moduli space of the reducible Spin$^c$ structure
$c$.
 This is certainly  expected.
Moreover,
 the following is proved in
\cite{FL} (see also similar statements in \cite{B}, \cite{BF} and
\cite {F}).

\begin{proposition} \label{relation} Let $M$ be a spin manifold with $2\chi+3\sigma=0$
and $b^+\geq 2$. Let  $c$ be a reducible Spin$^c$ structure.   Then,
$n_{M,c}=0$, and under the natural isomorphism between
$\Omega_0^{uo}$ and ${\mathbb Z}_2$, the unoriented bordism class is
equal to the ordinary SW invariant modulo $2$.
\end{proposition}

\begin{remark}Let $M$ be a spin manifold with $2\chi+3\sigma=0$
and $b^+\geq 2$, and let  $c$ be a reducible Spin$^c$ structure.
It follows from Proposition \ref{relation} that  the Mod 2 Seiberg-Witten invariant of $c$
depends only on $b^+$ and $\hbox{ind}\{D_a\}_J\in KQ(J)$.
\end{remark}

\section{Vanishing of  the unoriented bordism Seiberg-Witten
invariant}

In this section  we prove a vanishing result of the the unoriented
bordism Seiberg-Witten invariant.

Suppose $J=T^{4l-v}$
 with $l\geq 0$ and $0\leq v\leq 3$. Suppose $V_0$, $V_1$ are quaternionic
bundles over $J$ with $$\hbox{rank}_{\mathbb
C}V_0-\hbox{rank}_{\mathbb C}V_1=2p,$$ and $W_0$, $W_1$ are trivial
real bundles with
$$\hbox{rank}_{\mathbb R}W_0-\hbox{rank}_{\mathbb
R}W_1=-(4p+4l-v-1+\alpha)$$ for some integer $\alpha$.

\begin{proposition} \label{vanishing1} Let $J, V_0, V_1, W_0, W_1$ be as above. If
$p+l+\alpha>1$ and $p+l\geq 1$,
 then
$$e_1(V_0\oplus W_0, V_1\oplus W_1)=0.$$
\end{proposition}

\begin{proof} We first apply  (\ref{stability}),  the stability property of
$e_1$, to make the following reduction.

\begin{lemma} \label{model}  $e_1(V_0\oplus W_0, V_1\oplus W_1)$ is the same as
$$\left\{
\begin{array}{ll}
e_1(Q\oplus \underline {\mathbb H}^{p}\oplus \underline{\mathbb
R}^{v}, \underline{\mathbb H}^l\oplus \underline {\mathbb
R}^{4p+4l-1+\alpha}),
&\hbox{if $p\geq 0$,}\\
 e_1(Q\oplus \underline{\mathbb R}^{v},
\underline{\mathbb H}^{l-p}\oplus \underline {\mathbb
R}^{4p+4l-1+\alpha}),
&\hbox{if $p<0$,}\\
\end{array}\right.
$$
where  $Q$ is some quaternionic vector bundle with (complex) rank
$2l$.

\end{lemma}

\begin{proof}

By possibly stabilizing $V_0$ we can assume that $V_0$ has rank at
least $2l$. By Theorem \ref{monomorphism}, we can sum $V_1$ with a
quaternionic bundle to make it trivial. Therefore we can assume that
 $V_1= \underline {\mathbb H}^{a}$ for some positive integer $a\geq
-p$. Now, by Corollary \ref{split2}, if $p\geq 0$, we can write
$V_0=Q\oplus \underline {\mathbb H}^{a+p-l}$ where $Q$ is a complex
rank $2l$ quaternionic bundle. Similarly, if $p<0$, we can write
$V_0=Q\oplus \underline {\mathbb H}^{a+p-l}$.

\end{proof}
We assume now that $(V_0\oplus W_0,  V_1\oplus W_1)$ is of the form
as in Lemma \ref{model}.  The next step is to   choose a judicious
map to compute $e_1(V_0\oplus W_0, V_1\oplus W_1)$. Let us first
deal with the case that $v=0$. Notice that in this case the summand
$\underline {\mathbb R}^{v}$ in $V_0$ is trivial.

 Since $Q$ has rank $2l$, and
$$l + [\frac{4l+2}{4}]=2l,$$ by Theorem
\ref{monomorphism}, there exists a $Pin(2)-$equivariant monomorphism
$$ m=(m_1,..., m_{2l}):Q\longrightarrow  \underline{\mathbb H}^{2l}.$$
Write
$$\underline {\mathbb R}^{4p+4l-1+\alpha}
= \underline{\hbox{Im}{\mathbb H}}^{l+p}\oplus  \underline {\mathbb
R}^{l+p+\alpha-1}.$$ Consider the standard $Pin(2)-$equivariant
quadratic map
$$h: \underline{\mathbb H}\longrightarrow
\underline{\hbox{Im}{\mathbb H}}, \quad
      h(q)= qi\bar q.$$

When $p\geq 0$, we define, for $u\in Q$ and $(q_1,...,q_p)\in
\underline{\mathbb H}^p$, a $Pin(2)-$equivariant map
$$g_1:Q\oplus \underline{\mathbb H}^p\longrightarrow \underline{\mathbb H}^{2l}\oplus \underline{\mathbb H}^{p}
=\underline{\mathbb H}^l\oplus (\underline{\mathbb H}^{l}\oplus \underline{\mathbb H}^{p})
\longrightarrow \underline{\mathbb H}^l\oplus \underline{\hbox{Im}{\mathbb H}}^{l+p}$$
by
$$
g_1(u,q_1,...,q_p) =(m_1(u), ..., m_l(u),h(m_{l+1}(u)), ...,h(
m_{2l}(u)),h(q_1),..., h(q_p)).$$ And we define $g_0:Q\oplus
\underline{\mathbb H}^p\to  \underline{\mathbb H}^l\oplus
\underline{\hbox{Im}{\mathbb H}}^{l+p}$ to be the zero map. Then we
define, for $i=0,1$,  $$f_i:Q\oplus \underline{\mathbb H}^p
\longrightarrow \underline{\mathbb H}^l\oplus
\underline{\hbox{Im}{\mathbb H}}^{l+p}\oplus
 \underline {\mathbb R}^{l+p+\alpha-1}$$
by  $f_i=(g_i, k_i)$ with $k_1=(0,...,0)$ and $k_0=(1,...,1)$.
Notice that we can define $k_0$ this way since $l+p+\alpha$ is
assumed to be at least 2. Clearly $f_0$ is $U(1)-$equivariant and
non-vanishing. Since $g_1$ is $Pin(2)-$equivariant and non-vanishing
on the sphere bundle, so is $f_1$.   Moreover, we see that the
linear homotopy $\tilde f_t=tf_0+(1-t)f_1$  is never $0$ in the
$\underline {\mathbb R}^{l+p+\alpha-1}$ summand for $t\in (0,1)$.
Therefore $\tilde f_t^{-1}(0)=\emptyset$.

When $p\leq 0$,
we define, for $u\in Q$, a $Pin(2)-$equivariant map
$$g_1:Q\longrightarrow \underline{\mathbb H}^{2l}=\underline{\mathbb H}^{l-p}\oplus \underline{\mathbb H}^{l+p}
\longrightarrow \underline{\mathbb H}^{l-p}\oplus
\underline{\hbox{Im}{\mathbb H}}^{l+p}
$$
by
$$
g_1(u) =(m_1(u), ..., m_{l-p}(u),h(m_{l-p+1}(u)), h(m_{2l}(u))).$$
And we again define $g_0:Q\oplus \underline{\mathbb H}^p\to
\underline{\mathbb H}^l\oplus \underline{\hbox{Im}{\mathbb
H}}^{l+p}$ to be the zero map, and define in the same way, for
$i=0,1$,
$$f_i:Q \longrightarrow \underline{\mathbb H}^{l-p}\oplus
\underline{\hbox{Im}{\mathbb H}}^{l+p}\oplus
 \underline {\mathbb R}^{l+p+\alpha-1}$$
by  $f_i=(g_i, k_i)$ with $k_1=(0,...,0)$ and $k_0=(1,...,1)$. It is
easy to see that the linear homotopy $\tilde f_t$ has the same
property as in the case $p\geq 0$.

Notice that in both cases $f_1$ is
independent of $a\in J$ when restricted to $S(0\oplus W_0)$, i.e.
admissible. Thus we can use it to compute $e_1(V_0\oplus W_0,
V_1\oplus W_1)$. Since $\tilde f^{-1}_t(0)=\emptyset$ and $\tilde
f_t$ is  independent of $a\in J$ when restricted to $S(0\oplus
W_0)\times [0,1]$ we conclude that $e_1(V_0\oplus W_0, V_1\oplus
W_1)=0$ in the case of $v=0$.


 We now use the trick in \cite{L} to reduce the general case  to the case
 $v=0$.

Let $O_{T^v}$ be the point in $T^v$ coming from the origin of
${\mathbb R}^v$ and $B^v$ be  an invariant disc of $T^v$ around
$O_{T^v}$  and $\rho:B^v\longrightarrow {\mathbb R}^v$ be an
equivariant diffeomorphism. Consider the projection and the
embedding
$$p:T^{4l-v}\times T^v\longrightarrow
T^{4l-v},\quad  e:T^{4l-v}\longrightarrow T^{4l-v}\times O_{T^v}.
$$
Via $\rho$ we identify $ p^*Q|_{T^{4l-v}\times B^v}$ with the bundle
$Q\oplus \underline {\mathbb R}^v$   over $T^{4l-v}$.
 Notice that this identification is $Pin(2)-$equivariant since $\rho$ is.
Via this identification the monomorphism
$m:p^*Q\longrightarrow \underline{\mathbb H}^{2l} $ induces
a $Pin(2)-$equivariant bundle map (not a homomorphism)
$$m':Q\oplus \underline {\mathbb R}^v\longrightarrow \underline{\mathbb H}^{2l}$$
by the formula
$$m'|_z( u, s)=m|_{z\times \rho^{-1}(s)}(p^*u),$$
 where $z\in T^{b_1}, u\in Q$ and $s\in \underline{\mathbb R}^v$.

Now we define $g_1$ in the same way except replacing $Q$ by $Q\oplus
\underline {\mathbb R}^v$, replacing $m$ by $m'$ and adding a
monomorphism $\tau$ from ${\mathbb R}^v$ to the first
$\underline{\hbox{Im}\mathbb H}$.  We only need to verify that $g_1$
is $Pin(2)-$equivariant and non-vanishing on the sphere bundle,  the
remaining arguments are exactly the same as in the case $v=0$. $g_1$
is clearly $Pin(2)-$equivariant as  the linear map $\tau$ is
$Pin(2)-$equivariant. Since $m$ is a monomorphism the
$\underline{\mathbb H}^l$ component of $g_1$ vanishes only if $u=0$.
And if $u=0$ then the first $\underline{\hbox{Im}\mathbb H}$
component of $g_1$ vanishes only if $s=0$ as $h(0)=0$ and $\tau$ is
a monomorphism. Thus $g_1$ does not vanishes on the sphere bundle.
\end{proof}

\begin{theorem} \label{vanishing2}Let $M$ be a spin $4-$manifold with
\begin{equation}\label{notation}
b_1=4l-v, \,\,\, 0\leq v\leq 3, \quad \sigma=-16p, \quad
2\chi+3\sigma=4\alpha
\end{equation}
 and $b^+\geq 2$.  Let  $c$ be a reducible
Spin$^c$ structure.
 If $p+l+\alpha>1$ and $p+l\geq 1$, then
$SW^{uo}(M,c)$ is zero.
\end{theorem}
\begin{proof}
 By (\ref{euler}) and (\ref{2euler+3signature}), we have
$$-16p=\sigma=(b^{+}-b^{-})=-4(1-b_1+b^+)+4\alpha,$$ and hence
\begin{equation} \label{b+} b^+=4p+b_1-1+\alpha=4p+4l-1+\alpha-v.
\end{equation}
 Thus, any $(V_0, W_0, V_1,W_1)$ arising from an admissible finite
dimensional approximation of the Seiberg-Witten equations associated
to $c$ satisfies the assumption in  Proposition \ref{vanishing1}.
\end{proof}

\begin{remark}
In \cite{FL}, in the case $\alpha=0$ and $p+l=1$,  we are able to identify
 $SW^{uo}(M,c)$ with the $\epsilon$ invariant in
\cite{FKMM}.
\end{remark}

\section{Proof of  Theorem \ref{main}}

 Let us first recall  some relevant facts  (see \cite{L}, \cite{McS}, \cite{T})  about
minimal symplectic $4-$manifold with Kodaira dimension zero.

\begin{lemma} \label{fact} Let $(M,\omega)$  be a minimal symplectic $4-$manifold with
Kodaira dimension zero, then it has torsion canonical class
$K_{\omega}$. Moreover, it has the following properties.

\begin{enumerate}
\item  $2\chi+3\sigma=0$ and $M$ has even intersection form.

\item $K_{\omega}$ is either trivial, or  of order two which only
occurs when  $M$ is an integral homology Enriques surface. In
particular, $M$ is spin and the spin$^c$ structure ${\mathcal
K}_{\omega}$ is reducible except when $M$ is an integral homology
Enriques surface.

\item When $b^+\geq 2$, the Mod 2 Seiberg-Witten invariant of
${\mathcal K}_{\omega}$ is nonzero.

\end{enumerate}
\end{lemma}

We note that $b^-$ can be expressed via (\ref{2euler+3signature})
 as \begin{equation} \label{b-} b^-=4-4b_1+5b^+.
\end{equation}

Next we list minimal K\"ahler surfaces with $\kappa=0$ and
orientable $T^2-$bundles over $T^2$ in the following table according
to their homology type.

\begin{table}[ht]
\caption{}\label{eqtable}
\renewcommand\arraystretch{1.5}
\noindent\[
\begin{array}{|c|c|c|c|c|c|c|}
\hline class& b^+ & b_1 & \chi& \sigma & b^-
 & \hbox{known as }\\
\hline
 a)&3 & 0 &24&-16 & 19   & K3  \\
\hline
 b)& 3 & 4 &0&0 & 3 & \hbox{4-torus}  \\
\hline
 c) &2 & 3 &0& 0 & 2 &  \hbox{primary Kodaira
surface} \\
\hline
 d)&1 & 0 &12& -8 & 9 &\hbox{Enriques surface} \\
\hline e)&1 & 2 &0&0 & 1& \hbox{hyperelliptic
surface if complex}  \\
\hline
\end{array}
\]
\end{table}

We now finish the proof of Theorem \ref{main}.

\begin{proof} Let $M$ be a minimal symplectic
$4-$manifold with Kodaira dimension zero.

{\bf Bounds on $b^+, b^-$ and $b_1$}.
 Suppose $M$ is non-spin. In this case $M$ is
an integral homology Enriques surface by Lemma \ref{fact}. In
particular,
\begin{equation} \label{nonspin}
b^+=1, \quad b^-=9,\quad b_1=0.
\end{equation}

Since $M$ is symplectic, $b^+$ is at least $1$. Suppose $b^+=1$.
Then by (\ref{b-}) we have $b^-=9-4b_1$. Since $M$ has even
intersection form, we have $\sigma=b^+-b^-=1-b^-$ is divisible by
$8$. Moreover $b^-$ is non-negative, thus we have only two cases:
\begin{equation} \label{b+=1} b^+=1, \quad   b^-=9,\quad  b_1=0
\end{equation} or
\begin{equation} \label{b+=1'}
 b^+=1,\quad b^-=1,\quad b_1=2
 \end{equation}

Now we assume that $M$ is spin and has $b^+\geq 2$. We then can use
(\ref{notation}) for the homological invariants of  $M$. Notice that
  $\alpha=0$ by Lemma \ref{fact}.
 Then by the
vanishing from Theorem \ref{vanishing2}  and the non-vanishing from
Proposition \ref{relation} and Lemma \ref{fact}, we conclude that
\begin{equation} \label{bound} p+l\leq 1.
\end{equation}

On the other hand, we have by (\ref{b+})
 \begin{equation} \label{b+2}b^+=4p+4l-v-1.
 \end{equation}
Since $b^+$ is non-negative, we have
$$p+l\geq {1+v\over
4}.$$ Thus, as an integer, we must have $p+l\geq 1$.
 It then follows
from  (\ref{bound}) that
\begin{equation} \label{=} p+l=1.
\end{equation}

It follows from (\ref{b+2}) and (\ref{=}) that
\begin{equation} \label{b+3}
b^+=4-v-1\leq 3. \end{equation}

Since $b_1$ is non-negative,  we have from
 (\ref{b-}) and (\ref{b+3}) that
\begin{equation}\label{b-19}
b^-\leq 4+5b^+\leq 19.
\end{equation}
 Since  $b^-$ is non-negative, we have by
(\ref{b-}) and (\ref{b+3}) that
$$4b_1=4+5b^+-b^-\leq 4+5b^+\leq 19.$$  Hence we
conclude that
\begin{equation}\label{b14}
b_1\leq 4.
\end{equation}

Putting together  (\ref{nonspin}), (\ref{b+=1}), (\ref{b+=1'}),
(\ref{b+3}), (\ref{b-19}) and (\ref{b14}), we obtain the desired
Betti number bounds
$$b^+\leq 3, \quad b^-\leq 19, \quad b_1\leq 4.$$

{\bf Bounds on $\chi$ and $\sigma$}.  For the signature
$\sigma=b^+-b^-$, it is then between $-19$ and $3$. Since  $\sigma$
is divisible by $8$, $\sigma$ can only be $-16, -8, $ or $0$. It
follows from $2\chi+3\sigma=0$ that the Euler characteristic $\chi$
can only be $24, 12,$ or $0$.

{\bf Homology type}. Comparing with Table \ref{eqtable},
 in the case that $M$ is not spin or has $b^+=1$,   $M$ is either a homology
Enriques surface, or a homology $T^2-$bundle over $T^2$. In the case
that $M$ is spin and has $b^+\geq 2$, then $p+l=1$.  Since $l\geq 0$
and we have shown that $p=\frac{-\sigma}{16}\geq 0$, we have  either
$l=0, p=1$ or $l=1, p=0$. When $l=0, p=1$, $M$ is a homology K3.
 When $l=1, p=0$, we have $b_1\leq 4$ and $b^+=b^-=b_1-1$. In this case, $M$ is a homology
$T^2-$bundle over $T^2$ according to Table \ref{eqtable}.

\end{proof}

Finally the proof of  Corollary \ref{nonminimal}.

\begin{proof}
A non-minimal symplectic 4--manifold with $\kappa=0$ is obtained
from blowing up a minimal one. The blow up process keeps $b^+$ and
$b_1$ unchanged and increases $b^-$. Hence we still have the bounds
$b^+\leq 3$ and $b_1\leq 4$, as well as the bound $\chi\geq 0$.
 \end{proof}

\begin{remark}
In the broad context of the geography problem of symplectic
4-manifolds (see the survey \cite{L2}). Theorem \ref{main} and
Corollary \ref{nonminimal} provide complete answers in the case
$\kappa=0$.
\end{remark}

\end{document}